\documentclass[12pt,twoside,a4paper]{article}
\usepackage{amssymb,amsbsy,amsmath,amscd,eucal,exscale}

\newtheorem{theorem}{Theorem}[section]
\newtheorem{lem}[theorem]{Lemma}
\newtheorem{prop}[theorem]{Proposition}
\newtheorem{cor}[theorem]{Corollary}
 \title{Acyclicity of Schneider and Stuhler's coefficient systems: 
another approach in the level $0$ case}

\author{Paul Broussous\\
UMR 6086 du CNRS et  Universit\'e de Poitiers\\
 T\'el\'eport 2 - BP 30179\\
 Bd Pierre et Marie Curie\\
 86962 Futuroscope Chasseneuil\\
France\\
Email: broussou@mathlabo.univ-poitiers.fr}

\begin{document}
\maketitle

\begin{abstract}

Let $F$ be a non archimedean local field and $G$ be the locally
profinite group ${\rm GL}(N,F)$, $N\geq 1$. We denote by $X$ the
Bruhat-Tits building of $G$. For all smooth complex representation
${\mathcal V}$ of $G$ and for all level $n\geqslant  1$, 
Schneider and Stuhler  have constructed a coefficient system ${\mathcal C}
= {\mathcal C}({\mathcal V} ,n)$ on the simplicial
complex $X$. They proved that if $\mathcal V$ is generated by its fixed
vectors under the principal congruence subgroup of level $n$, then the
augmented complex $C_{\bullet}^{\rm or}(X,{\mathcal C} )\rightarrow
{\mathcal V}$ of oriented chains of $X$ with coefficients in $\mathcal
C$ is a resolution of $\mathcal V$ in the category of smooth complex
representations of $G$. In this paper we give another proof of
this result, in the level $0$ case, and assuming moreover that
$\mathcal V$ is generated by its fixed vectors under an Iwahori subgroup $I$ of
$G$.  Here ``{\it level} $0$'' refers to Bushnell and Kutzko's
terminology, that is to the case $n=1+0$.  Our approach is different.
We strongly use the fact that the trivial character of $I$ is a type
in the sense of Bushnell and Kutzko.
\end{abstract}



\section{Introduction} 

Let $F$ be a locally compact non-archimedean  field and $G$
be the group ${\rm GL}(N,F)$, where $N\geqslant 1$ is a fixed integer. We
denote by $X$ the semisimple affine building of $G$. We let ${\mathcal
R}(G)$ denote the category of smooth complex representations of
$G$. For any simplex $\sigma$ of $X$, we denote by $U_{\sigma}$ the
pro-unipotent radical of the (compact) parahoric subgroup $P_{\sigma}$
attached to $\sigma$. We fix a representation ${\mathcal V}\in {\mathcal R}(G)$ and
we assume that it is generated by its fixed vector under some Iwahori
subgroup $I$. In \cite{[SS]} Schneider and Stuhler define a coefficient
system ${\mathcal C} =({\mathcal V}_{\sigma})_{\sigma}$, where
${\mathcal V}_{\sigma}={\mathcal V}^{U_{\sigma}}$ for all simplex $\sigma$ of $X$. The
central result of {\it loc. cit.} is that the augmented complex of
oriented chains of $X$ with coefficients in ${\mathcal C}$

\begin{equation}\label{a}
0\rightarrow C_{N-1}^{\rm or}(X,{\mathcal C} ) {\buildrel \partial \over
\rightarrow} \cdots{\buildrel \partial \over \rightarrow} 
C_{0}^{\rm or}(X,{\mathcal C} ){\buildrel \epsilon \over \rightarrow} {\mathcal V}
\end{equation}

is a resolution of ${\mathcal V}$ in ${\mathcal R}(G)$. Of course the result of
\cite{[SS]} is more general and is  concerned with  representations of any level. 

The aim of this paper is to give another proof of the acyclicity of
the homological complex (1). Our proof proceeds as follows. Let ${\mathcal
R}_{(I,1)}(G)$ be the full subcategory of ${\mathcal R}(G)$ whose
representations are generated by their fixed vectors under $I$. Let
${\mathcal H}(I,1)$ be the Hecke algebra of the trivial character of
$I$. Then since $(I,1)$ is a type for $G$ (see \cite{[BK1]}, \cite{[BK2]}), we have
an equivalence of categories
$$
{\mathcal R}_{(I,1)}(G) \longrightarrow {\rm Mod}({\mathcal H}(I,1))
$$
$$
W\mapsto W^{I}
$$
where ${\rm Mod}({\mathcal H}(I,1))$ is the categoy of left unital modules
over ${\mathcal H}(I,1)$. We then show that the complex (\ref{a}) is actually a
complex in ${\mathcal R}_{(I,1)}(G)$. Applying the equivalence of categories
to (\ref{a}), we show that we obtain a complex isomorphic to the chain
complex of the standard apartment of $X$ with constant coefficients in
${\mathcal V}^I$. Since this apartment is contractible, this latter complex is
exact.
\smallskip

 The paper is organized as follows. In {\S}7 we compute the image of
the complex (1) under the equivalence of categories. In {\S}8 we show
how this image is related to the chain complex of the standard
apartment with coefficients in ${\mathcal V}^I$. Proofs are based on
some geometrical lemmas (cf.  {\S}5) and on a key lemma on
idempotents of the Hecke algebra of $G$ (cf. {\S}6).

I thank Guy Henniart and Peter Schneider for stimulating
discussions. The proof of the main result of {\S}6 is due to
Jean-Fran\c cois Dat.

\section{ Notation} We keep the notation of the introduction.
  We shall denote by $Y_q$
the set of $q$-dimensional simplices of a simplicial complex $Y$,
$q\in {\mathbb Z}_{\geqslant 0}$. We  set
$$
G^o =\{g\in G\ ; \ {\rm det}(g)\in {\mathfrak o}_{F}^{\times}\}\ ,
$$
where ${\mathfrak o}_F$ denotes the ring of integers on $F$. 
We fix a maximal split torus $T$ of $G$ and denote by $\mathcal A$ the
corresponding appartement of $X$. We fix a chamber $C_o$ of $\mathcal
A$ and write 
$$
I=\{ g\in G^o\ ; \ gC_o =C_o\}
$$
for the Iwahori subgroup attached to $C_o$. The subgroup $I$ together
with
$$
N^{o}=\{g\in G^o \ ; \ gTg^{-1}=T\}
$$
form an affine $BN$-pair in $G^o$. The corresponding affine Weyl group
is $W^{\rm aff}=N^o / T^o$, where $T^o \subset T$ is the maximal
compact subgroup of $T$. Recall that $X$ is the geometric realization
of the combinatorial building of the $BN$-pair $(I, N^o )$.

\smallskip

 We fix a labelling $\lambda$~: $X\rightarrow \Delta_{N-1}$ (cf. \cite{[Br]},
 page 29), where
 $\Delta_{N-1}$ is the standard $(N-1)$-simplex. The action of $G^o$
 is label preserving. The labelling $\lambda$ gives rise to an
 orientation of the simplices of $X$ and to incidence numbers $[\sigma
 :\tau ]$ , $\tau\subset \sigma$, ${\rm dim}\sigma ={\rm dim}\tau +
 1$. The action of $G^o$ preserves the orientation and the incidence
 numbers.
\smallskip

 Finally we fix a smooth representation ${\mathcal V}$ in the subcategory
 ${\mathcal R}_{(I,1)}(G)$ of ${\mathcal R}(G)$. We have a $G$-equivariant
 coefficient system 
$$
{\mathcal C} =({\mathcal V}_{\sigma})_{\sigma}
$$
on $X$, where ${\mathcal V}_{\sigma}={\mathcal V}^{U_{\sigma }}$, for all
simplex $\sigma$. Here the restriction maps $r_{\tau}^{\sigma}$~:
${\mathcal V}_{\sigma}\rightarrow {\mathcal V}_{\tau}$, for $\tau\subset \sigma$, are
the inclusion maps.
\bigskip

\section{The chain complexes}
\smallskip

 The coefficient system ${\mathcal C}$ gives rise to the augmented  complex of
 oriented chains of $X$ with coefficients in ${\mathcal C}$: 

\begin{equation}\label{b}
0\rightarrow C_{N-1}^{\rm or}(X,{\mathcal C} ){\buildrel \partial \over \rightarrow}  \cdots {\buildrel \partial \over \rightarrow}
C_{0}^{\rm or} (X, {\mathcal C}
){\buildrel \epsilon  \over \rightarrow}{\mathcal V}
\end{equation}

We do not repeat the definition and we refer to \cite{[SS]}{\S}3, for we
shall actually use another chain complex.

For $q=0,\dots ,N-1$, we consider the $G$-module
$$
C_{q}(X,{\mathcal C}) =\bigoplus_{{\rm dim}\sigma =q}{\mathcal V}_{\sigma}
$$
where the action of $g\in G$ is given by $g.(v_{\sigma})_{\sigma}
=(gv_{g^{-1}\sigma})_{\sigma}$. We have boundary maps
$$
\partial~: C_{q+1}(X, {\mathcal C} )\rightarrow C_{q}(X,{\mathcal C} )\ , \ (v_{\tau}
)_{\tau}\mapsto (w_{\sigma})_{\sigma}\ ,
$$
where 
$$
w_{\sigma}=\sum_{\tau\supset \sigma \ , \ {\rm dim}\tau =q+1}[\tau
:\sigma ]v_{\tau}\ .
$$
These boundary maps are $G^o$-equivariant. We have an augmented complex
of chains in the category of $G^o$-modules:

\begin{equation}\label{c}
0\rightarrow C_{N-1}(X,{\mathcal C} ){\buildrel \partial \over
\rightarrow}  \cdots {\buildrel \partial \over \rightarrow}  C_{0}(X, {\mathcal C}
){\buildrel \epsilon \over \rightarrow}{\mathcal V}
\end{equation}

where $\displaystyle \epsilon [(v_{\sigma})_{\sigma}]=\sum_{{\rm
dim}\sigma =0} v_{\sigma}$. 

 This complex is isomorphic to the augmented complex of oriented chains
(\ref{b})   when we see the latter as a complex
of $G^o$-modules.
\smallskip

\begin{prop}
 The chain complex (\ref{b}) is a chain complex in
the category ${\mathcal R}_{(I,1)}(G)$.
\end{prop}

{\bf Proof}. Using the labelling $\lambda$, any simplex $\sigma$ in
$X$ gives rise to two oriented
simplices (cf. \cite{[SS]}{\S}3)  $\langle \sigma \rangle$ and $\langle {\bar \sigma}\rangle$
(where $\langle \sigma \rangle = \langle {\bar \sigma}\rangle$ if
$\sigma$ is a vertex). For all $q\in \{ 0,\dots ,N-1\}$ and for all
$q$-simplex $\sigma$, we write
$C_{q}^{\rm or}(X ,\sigma )$ for the space of oriented $q$-chains with
support in $\{ \langle \sigma \rangle\cup \langle {\bar
\sigma}\rangle\}$. Then we have a natural isomorphism of
$P_{\sigma}$-modules $C_{q}^{\rm or}(X ,\sigma )\simeq {\mathcal V}_{\sigma}$.

Fix $q\in \{ 0,\dots ,N-1\}$. We must prove that, as a
$G$-module, $C_{q}^{\rm or}(X,{\mathcal C} )$ is generated by a subspace fixed by $I$. Fix a vertex
$s$ of $C_o$. Then any $q$-simplex of $G$ has a $G$-conjugate $\sigma$
satisfying $s\subset \sigma \subset C_o$. We fix a system $\Sigma$ of
representatives of the $G$-conjugacy classes of $q$-simplices $\sigma$
satisfying $s\subset \sigma \subset C_o$. We have the decomposition
$$
C_{q}^{\rm or}(X ,C)=\bigoplus_{\sigma\in \Sigma}\bigoplus_{g\in
G/G_{\sigma}}gC_{q}^{\rm or}(X ,\sigma ) \ ,
$$
where $G_{\sigma}\supset P_{\sigma}$ is the stabilizer of $\sigma$ in
$G$. So we are reduced to proving that for all $\sigma\in \Sigma$, the
$G$-module $\displaystyle \bigoplus_{g\in G/G_{\sigma}}gC_{q}^{\rm
or}(X ,\sigma ) $ lies in ${\mathcal R}_{(I,1)}(G)$. For this it suffices
to prove that $C_{q}^{\rm or}(X ,\sigma )$  is generated by
$C_{q}^{\rm or}(X ,\sigma )^{I}$ as a $P_{\sigma}$-module, in other
words that ${\mathcal V}_{\sigma}$ is generated by ${\mathcal V}_{\sigma}^{I}$ as a
$P_{\sigma}$-module.  

We may see ${\mathcal V}_s :={\mathcal V}^{U_{s}}$ as a complex representation of
${\mathbb G} =P_s /U_s \simeq {\rm GL}(N,{\bf k})$, where $\bf k$ is the residue
field of $F$. Moreover the image  ${\mathbb P}_{\sigma}$
(resp. $\mathbb B$) of
$P_{\sigma}$ (resp. $I$) in $P_s /U_s$ is a parabolic subgroup of
$\mathbb G$ (resp. a Borel subgroup). The unipotent radical
${\mathbb U}_{\sigma}$ of ${\mathbb P}_{\sigma}$ is the image of $U_{\sigma}$ in
$P_s /U_s$. In particular we may see
${\mathcal V}_{\sigma}=({\mathcal V}_{s})^{U_{\sigma}}$ as the Jacquet module of
${\mathcal V}_{s}$ with respect to ${\mathbb U}_{\sigma}$. We have the following
key result due to Schneider and Zink:
\smallskip

\begin{lem} (\cite{[SZ]} Prop. (5.3), page 19.) The
$\mathbb G$-module ${\mathcal V}_{s}$ has cuspidal support $({\mathbb T} ,1)$, where
$\mathbb T$ is the Levi subgroup of ${\mathbb B}$ (a maximal split torus in
$\mathbb G$).
\end{lem}

It follows that ${\mathcal V}_{\sigma}={\mathcal V}_{s}^{\mathbb U}$ has cuspidal support
$({\mathbb T} ,1)$ as well. By Frobenius reciprocity for parabolic
induction, we have that ${\mathcal V}_{\sigma}$ is generated by
${\mathcal V}_{\sigma}^{\mathbb B}=({\mathcal V}_{\sigma})^{I}$ as a ${\mathbb P} (=   
P_{\sigma}/U_{\sigma})$-module and we are done.

\section{The strategy}

We know that (\ref{b}) is a chain complex in ${\mathcal R}_{(I,1)}(G)$. Let
${\mathcal H}(I,1)$ denote the Hecke algebra of $(I,1)$ and ${\rm
Mod}({\mathcal H}(I,1))$ the category of left unital ${\mathcal
H}(I,1)$-modules. We have an equivalence of categories
(cf. \cite{[BK1]}, \cite{[BK2]}):
$$
{\mathcal R}_{(I,1)}(G)\ {\buildrel \Phi\over\longrightarrow}\ {\rm
Mod}({\mathcal H}(I,1))
$$
$$
{\mathcal W}\mapsto {\mathcal W}^{I} \ .
$$
So in order to prove that the augmented chain complex (\ref{b}) is acyclic it
suffices to prove that the corresponding augmented chain complex in 
${\rm Mod}({\mathcal H}(I,1))$ is acyclic. We must prove that 
$$
0\rightarrow C_{N-1}^{\rm or}(X,{\mathcal C} )^{I}\ {\buildrel \Phi (\partial
)\over\longrightarrow}\  \cdots\  {\buildrel \Phi (\partial
)\over\longrightarrow}\  C_{0}^{\rm or}(X,{\mathcal C} )^{I}\  
 {\buildrel \Phi (\epsilon )\over\longrightarrow} \ {\mathcal V}^{I}
$$
is exact. But since $I\subset G^o$ and the complexes (\ref{b}) and (\ref{c})
are isomorphic as complexes of $G^o$-modules, we are reduced to proving:
\smallskip

\begin{prop} The following sequence of $\mathbb C$-vector
spaces is exact:
\begin{equation}\label{d}
0\rightarrow C_{N-1}(X,{\mathcal C} )^{I}\ {\buildrel \Phi (\partial
)\over\longrightarrow}\  \cdots\  {\buildrel \Phi (\partial
)\over\longrightarrow}\  C_{0}(X,{\mathcal C} )^{I}\  
 {\buildrel \Phi (\epsilon )\over\longrightarrow} \ {\mathcal V}^{I}\ .
\end{equation}

\end{prop}

The proof will proceed in several steps.

\section{Some geometric lemmas}

 We fix a simplex $\sigma\subset C_o$ and an element $w\in W^{\rm
 aff}$. Let $E[\sigma ,w^{-1}C_o ]$ be the {\it enclos} of $\sigma\cup
 w^{-1}C_o$ (cf. \cite{[BT]}, Definition (2.4.1)).
 This is a convex subset of $\mathcal A$ as well as a
 subsimplicial complex.
\smallskip

\begin{lem} There exists a unique chamber $C_1 = C_1
(\sigma ,w)$ of $\mathcal A$ such that
$$
\sigma \subset C_1 \subset E[\sigma ,w^{-1}C_o ]\ .
$$
\end{lem}

{\it Proof}. The element $w$ does not play any special role here and
we may assume that $C_2 =w^{-1}C_o$ is any chamber of $\mathcal A$. Recall
\cite{[BT]} that $E[\sigma ,C_2 ]$ is a union of closed chambers. This
proves the existence of $C_1$. By \cite{[BT]}, Prop (2.4.5),
$E[\sigma ,C_2 ]$ is the intersection of the half-apartments
(corresponding to affine roots) containing both $\sigma$ and
$C_2$. Assume that $C_1$ and $C_1 '$ are two distinct chambers satisfying the
assumption of the lemma. Since the set $E[\sigma ,C_2 ]$ is closed in
the sense of \cite{[BT]}, it contains all chambers in a minimal
gallery connecting $C_1$ and $C_1 '$. We may choose such a gallery so
that its chambers contain $\sigma$. It follows that we may assume
that $C_1$ and $C_1 '$ are adjacent (and contain $\sigma$). Let $H$ be
a  half-apartment whose hyperplan contains $C_1 \cap C_1 '\supset
\sigma$. We may choose $H$ so that $H$ contains $C_2$. But then $H$
must contain the enclos $E[\sigma ,C_2 ]\supset C_1 \cup C_2$, a
contradiction.

\smallskip

We let $C(\sigma ,w )$ denote the chamber $wC_{1}$. By equivariance
and using the last lemma, $C(\sigma ,w)$ is the unique chamber of
$\mathcal A$ satisfying 
$$
w\sigma \subset C(\sigma ,w )\subset E[w\sigma ,C_o ] \ .
$$
 It only depends on the simplex $w\sigma$ (and on the
reference chamber $C_o$).

\begin{lem} With the previous notation, we may choose
points $x_{\sigma}\in \vert\sigma\vert^o$, $x_1\in \vert C_1\vert^o$,
$x_{w}\in \vert w^{-1}C_o\vert^o$, in the geometric realizations, such
that $x_1$ belongs to the geodesic segment $[x_{\sigma}x_{w}]$.
\end{lem}

{\it Remark}. In the terminology  of \cite{[SS]}{\S}2, this means that the
simplex $C_1$ lies between $\sigma$ and $w^{-1}C_o$.
\smallskip

{\it Proof}. Let $x_{\sigma}$ be the isobarycenter of $\sigma$ and
${\mathfrak C}$ be the union of the open segments $(xx_{\sigma})$ where $x$ runs
over $\vert w^{-1}C_o\vert^o$. If $x_{\sigma}$ is the isobarycenter
of $w^{-1}C_o$, then $\sigma =C_1 =w^{-1}C_o$ and the lemma is
clear. So we assume that $x_{\sigma}$ is not the isobarycenter of
$w^{-1}C_o$ so that  ${\mathfrak C}$ is an open subset of $\mathcal
A$, contained in the enclos $E[\sigma ,w^{-1}C_o ]$. There is an open
chamber $\vert D\vert^o$,   contained in  $E[\sigma ,w^{-1}C_o ]$,
intersecting ${\mathfrak C}$ and such that $x_{\sigma}\in D$; otherwise  we
would have a covering of ${\mathfrak C}$ by  chambers $ D_{i}$,
$î=1,\dots ,r$, such that 
$$
x_{\sigma}\not\in \bigcup_{i=1,\dots ,r}D_{i}\supset {\mathfrak C}
$$
and this would contradict the fact that $x_{\sigma}$ is in the closure
of ${\mathfrak C}$. But since $x_{\sigma}\in D$, we have $\sigma\subset D$ and
by unicity in lemma 4, we have $D=C_1$.

\begin{lem} The subgroup $\langle U_{\sigma} ,
P_{\sigma}\cap I^{w^{-1}}\rangle$ generated by
$U_{\sigma}$ and $P_{\sigma}\cap I^{w^{-1}}$ is the
Iwahori subgroup $I_1$ wich fixes $C_1$.
\end{lem}

{\it Proof}.  We have $U_{\sigma}\subset  I_1$,
 since $\sigma\subset C_1$. Moreover $P_{\sigma}\cap I^{w^{-1}}$
 fixes the enclos $E[\sigma ,w^{-1}C_o ]$ pointwise and therefore
 fixes the chamber $C_1$. Hence $P_{\sigma}\cap I^{w^{-1}}
 \subset I_1$ and we have 
$$
\langle U_{\sigma} , P_{\sigma}\cap I^{w^{-1}}\rangle
 \subset I_1 \ .
$$
Let $I_{1}^{1}$ (resp. $I^1$) denote the pro-unipotent radical of $I_1$
(resp. of $I$). Since $C_1$ lies between $\sigma$ and $w^{-1}C_o$, 
by \cite{[SS]} Prop. (2.5), we have
$$
{\rm Im}[I_{1}^{1}\cap P_{\sigma}\rightarrow
P_{\sigma}/U_{\sigma}]
\subset
{\rm Im}[(I^{1})^{w^{-1}}\cap P_{\sigma}\rightarrow
 P_{\sigma}/U_{\sigma}]
$$
We therefore have the containment:
$$
U_{\sigma}[(I^{1})^{w^{-1}}\cap
P_{\sigma}]U_{\sigma}
\supset
U_{\sigma}[(I_{1}^{1})\cap
P_{\sigma}]U_{\sigma}\ .
$$
Since $I=T^o I^1$ and $I_1 =T^o I_{1}^{1}$, we get
$$
U_{\sigma}[I^{w^{-1}}\cap
P_{\sigma}]U_{\sigma}
\supset U_{\sigma}[ I_1 \cap P_{\sigma}]
U_{\sigma} =I_1\ .
$$
So we have that $I_1$ is contained in  $\langle U_{\sigma} ,
P_{\sigma}\cap I^{w^{-1}}\rangle$, as required.
\smallskip

 We shall need the following lemma.

\begin{lem} Let $K$, $J$, $L$ be Iwahori subgroups of $G$ whose
respective chambers $C_K$, $C_J$, $C_L$ lie in the apartment $\mathcal A$.
Assume that $C_J$ lies in the enclos $E[C_K ,C_L ]$ of $C_K\cup
C_L$. Then we have $J=(J\cap K)(J\cap L)$. In particular we have
$KJL=KL$.
\end{lem}

{\it Proof}. Let $\Phi$ (resp. $\Phi^{\rm aff}$ denote the set of roots
(resp. of affine roots) of $G$ relative to $T$. We have the
standard  valuated root datum $(T, (U_{{\bf a}})_{{\bf a}\in \Phi})$ of
$G$ (cf. \cite{[BT]}, {\S}6  and {\S}(10.2)). In particular for
all $a\in \Phi^{\rm aff}$, we have a congruence subgroup $U_{a}\subset
U_{\bf a}$, where $\bf a$ is the vector part of $a$. For any subset
$\Omega$ of $\mathcal A$, we denote by $\Phi^{\rm aff}(\Omega )$ the set
of affine roots which are positive on $\Omega$ and minimal for this
property, i.e. $a\in \Phi^{\rm aff}(\Omega )$ if $a(\Omega )\subset 
{\mathbb R}_{+}$ and $a(\Omega )-1\not\subset {\mathbb R}_{+}$. From {\it
loc. cit.} {\S}6,7, we have that  $J$ is generated by $T^o$ and the $U_a$,
where $a$ runs over $\Phi^{\rm aff}(C_J )$. More precisely the product
map
$$
T^o \prod_{a\in \Phi^{\rm aff}(C_J )}U_{a}\longrightarrow J
$$
is surjective for any ordering of the affine roots in $\Phi^{\rm aff}(C_J )$.
 On the other hand since $J\cap K$
(resp. $J\cap L$) is the fixator of the enclos $E[C_J, C_K ]$
(resp. $E[C_J ,C_L ]$) in $G^o$, we have that $J\cap K$ (resp. $J\cap
L$) is generated by $T^o$ and the $U_a$, $a$ running over $\Phi^{\rm
aff}(E[C_J, C_K ])$ (resp. $\Phi^{\rm aff}(E[C_J, C_L ])$). We claim
that any $a\in \Phi^{\rm aff}(C_J )$ is positive on $C_K$ or on $C_L$
(whence on $E[C_J , C_K]$ or on $E[C_J ,C_L ]$);
this will give $J\subset (J\cap K)(J\cap L)$, as required. Assume that
there exists $a\in \Phi^{\rm aff}(C_J )$ such that $a$ is negative on
$C_K$ and $C_L$. Then there would exist a  half-apartment with hyperplan ${\rm
ker}({\bf a})$ which contains  $C_K$ and $C_L$  but not $C_J$. But
this would contradict the fact that $C_J$ lies in the enclos of
$C_K\cup C_L$.

\section{A lemma on idempotents}

Let ${\mathcal H}(G)$ be the convolution algebra of locally constant
complex functions on $G$ with compact support, that is the {\it Hecke
algebra} of $G$. We shall denote by $\star$ the convolution product.
 For any Iwahori subgroup $K$ of $G$, we let $e_K$
denote the idempotent of ${\mathcal H}(G)$
attached to the trivial character of $K$. 
\smallskip

\begin{lem} Let $\mathcal W$ be a smooth representation of
$G$ and  $J$ be an  Iwahori subgroup of
$G$ whose chamber $C_J$ lies in $\mathcal A$. Then the map
${\mathcal W}^J \rightarrow {\mathcal W}^I$, $v\mapsto e_I .v$ is an
isomorphism of vector spaces.
\end{lem}

\begin{lem} With the notation as above, let $K$ be an
Iwahori subgroup of $G$ whose chamber $C_K$ lies in the enclos $E[C_J
,C_I ]$. Then the map ${\mathcal W}^J \rightarrow {\mathcal W}^I$, $v\mapsto
e_I .v$ factors through ${\mathcal W}^J\rightarrow {\mathcal W}^K\rightarrow
{\mathcal W}^I$, where the first map is given by $v\mapsto e_K .v$ and
the second is given by $w\mapsto e_I .w$.
\end{lem}

{\it Proof}. We claim that $e_I\star e_K\star e_J =e_I\star e_J$. The lemma will
follow easily. The support of $e_I\star e_K\star e_J$ is $IKJ$ wich is
$IJ = {\rm supp}(e_I\star e_J )$ by lemma 7. It suffices to prove
that $ e_I\star e_K\star e_J (1_G ) = e_I\star e_J (1_G )$, a
straightforward computation.

\smallskip

 As a consequence of this lemma, we have that in lemma 8, we may
 reduce to the case where the chambers $C_I$ and $C_J$ are adjacent.

\begin{lem}  With the notation of lemma 8, assume
moreover that the chambers $C_I$ and $C_J$ are adjacent. Then there
exist functions $f,g\in {\mathcal H}(G)$ such that:
$$
e_J =f\star e_I\star e_J\ {\rm and}\ e_{I}=e_{I}\star e_J\star g\ .
$$
\end{lem}

Let us first show how this latter lemma implies lemma 8 in the case
where $C_I$ and $C_J$ are adjacent. Let $v\in {\mathcal W}^J$. If
$e_I. v=0$,  we have $ f\star e_I\star e_J .v=e_J. v=0$, whence $v=0$ and
$f$ in injective. Let $w\in {\mathcal W}^I$. We have $w=e_{I}.w=e_{I} (e_J\star
g.w)$, with $e_J\star g.w\in {\mathcal W}^J$ and $e_I$ is therefore
surjective. 
\smallskip

{\it Proof of lemma 10}. (Compare \cite{[HL]} Theorem (2.5)). If
$N=1$, the result is trivial and we assume $N\geq 2$. We give the
proof of the existence of $f$; the proof of the existence of $g$ being
similar. Set $\tau =C_I \cap C_J$. This is a simplex of codimension
$1$. We have $I,J\subset P_{\tau}$ so that we may see $e_I$ and $e_J$ 
as idempotents in the group algebra of $P_{\tau}/U_{\tau}$. We fix an
isomorphism $P_{\tau}/U_{\tau}\simeq {\rm GL}(2,{\bf k})\times{\mathbb
T}$,
where ${\mathbb T}\simeq ({\bf k}^{\times})^{N-2}$, so that the image of $I$
(resp. $J$) in  $P_{\tau}/U_{\tau}$ is ${\mathbb B}^+ \times{\mathbb T}$
(resp. ${\mathbb B}^- \times {\mathbb T}$), where ${\mathbb B}^+$ (resp. ${\mathbb B}^-$) is
the Borel subgroup of upper (resp. lower) triangular matrices. So we
identify $e_I \simeq e_{{\mathbb B}^+\times {\mathbb T}}$, $e_J \simeq
e_{{\mathbb B}^-\times {\mathbb T}}$. Write ${\mathbb U}^+$ and ${\mathbb U}^-$ for the
unipotent radicals of ${\mathbb B}^+$ and ${\mathbb B}^-$ respectively. We have 
$$
e_{{\mathbb B}^+\times {\mathbb T}}\star e_{{\mathbb B}^-\times {\mathbb T}}\star
e_{{\mathbb B}^+\times {\mathbb T}}
=1/\vert {\mathbb B}^- \times {\mathbb T}\vert\sum_{x\in
{\mathbb B}^-\times{\mathbb T}}e_{{\mathbb B}^+\times {\mathbb T}}\star x  e_{{\mathbb B}^+\times
{\mathbb T}}
$$
$$
=1/q(q-1)^N\sum_{u\in {\mathbb U}^-}e_{{\mathbb B}^+\times {\mathbb T}}\star u e_{{\mathbb B}^+\times
{\mathbb T}}
$$
$$
=1/q(q-1)^N ( e_{{\mathbb B}^+\times {\mathbb T}} +  \sum_{u\in {\mathbb U}^-
, \ u\not= 1 }e_{{\mathbb B}^+\times {\mathbb T}}\star u e_{{\mathbb B}^+\times
{\mathbb T}} )\ .
$$
Here $q=\vert {\bf k}\vert$ is the size of the residue field. Any
$u\in {\mathbb U}^-\backslash \{ 1\}$ writes $u_{+}tsu_{+}'$, where
$u_{+},u_{+}'\in {\mathbb U}^+$, $t$ belongs to the diagonal torus of ${\rm
GL}(2,{\bf k})$ and $s$ is the matrix
$$
s=\left( 
\begin{matrix}
 0 & 1 \\
  -1 & 0
\end{matrix}
\right)\ .
$$
So we have 
$$
e_{{\mathbb B}^+\times {\mathbb T}}\star e_{{\mathbb B}^-\times {\mathbb T}}\star
e_{{\mathbb B}^+\times {\mathbb T}} = 
1/q(q-1)^N ( e_{{\mathbb B}^+\times {\mathbb T}} +  \sum_{u\in {\mathbb U}^-
, \ u\not= 1 }e_{{\mathbb B}^+\times {\mathbb T}}\star s e_{{\mathbb B}^+\times
{\mathbb T}} )\ .
$$
But since $e_{{\mathbb B}^+ \times {\mathbb T}}s=se_{{\mathbb B}^-\times{\mathbb T}}$, we get
$$
e_{{\mathbb B}^+\times {\mathbb T}}\star e_{{\mathbb B}^-\times {\mathbb T}}\star
e_{{\mathbb B}^+\times {\mathbb T}} = 
1/q(q-1)^N ( e_{{\mathbb B}^+\times {\mathbb T}} +
(q-1)se_{{\mathbb B}^-\times{\mathbb T}}\star e_{{\mathbb B}^+\times {\mathbb T}} )\ ,
$$
that is 
$$
(q(q-1)^Ne_{{\mathbb B}^+\times {\mathbb T}}- (q-1)s {\bf 1}_{{\rm GL}(2,{\bf
k})\times {\mathbb T}})\star  e_{{\mathbb B}^-\times {\mathbb T}}\star
e_{{\mathbb B}^+\times {\mathbb T}} = e_{{\mathbb B}^+\times {\mathbb T}}\ .
$$
In terms of functions on $G$, we have proved:
$$
(q(q-1)^N e_I -(q-1)e_{U_{\tau}sU_{\tau}})\star e_{J}\star e_I =e_I 
$$
and we are done.

\section{Computation of the $I$-fixed vectors}

Fix $q\in \{ 0,\dots ,N-1\}$. A set of representatives of the
$G^o$-conjugacy classes of $q$-simplices in $X$ is given by
$$
\Sigma_q :=\{\sigma_o \ {\rm simplex \ of}\ X; \ \sigma_o \subset C_o ,
\ {\rm dim}\sigma_o =q\} \ .
$$
If a simplex $\sigma$ is $G^o$-conjugate to $\sigma_o\in \Sigma_q$, we
say that $\sigma$ has type $\sigma_o$. The set of simplices of types
$\sigma_o$ is isomorphic to $G^o /P_{\sigma_o}$ as a
$G^o$-set. Any simplex of type $\sigma_o$ has a unique $I$-conjugate
lying in the apartment $\mathcal A$, so that we have the identification:
$$
I\backslash G/P_{\sigma_o}\simeq \{ \sigma\subset {\mathcal A}\ ; \
\sigma\ {\rm has \ type}\ \sigma_o \}\ .
$$
On the other hand, $W^{\rm aff}$ acts transitively on the simplices of
types $\sigma_o$ in  $\mathcal A$ so that we have: 
$$
W^{\rm aff}/W^{\rm aff}_{\sigma_o}\simeq \{\sigma\subset {\mathcal A}\ ; \ 
\sigma \ {\rm has  \ type}\ \sigma_{o}\} \ ,   
$$
where $W^{\rm aff}_{\sigma_o}=\{w\in W^{\rm aff}\ ; \ w\sigma_o
=\sigma_o\}$. Of course these identifications are compatible with the
canonical bijection
$$
W^{\rm aff}/W^{\rm aff}_{\sigma_o} \rightarrow I\backslash
G/P_{\sigma_o}
$$
$$
w\ {\rm mod} \ W^{\rm aff}_{\sigma_o}\mapsto Iw P_{\sigma_o}\ .
$$
For each $\sigma_o\in \Sigma_q$, we let $R_{\sigma_o}$ denote a set of
representatives of $W^{\rm aff}/W^{\rm aff}_{\sigma_o}$. We see
$W^{\rm aff}$ as a subgroup of $N^o$ in the usual way, so that
$R_{\sigma_o}\subset N^o$.

We have the decomposition
$$
C_{q}(X, {\mathcal C} )=\bigoplus_{{\rm dim}\sigma =q}{\mathcal V}_{\sigma}
=\bigoplus_{\sigma\in \Sigma_{q}} \bigoplus_{g\in G^o
/P_{\sigma}} {\mathcal V}_{g\sigma}
$$
and, since for $\sigma\in \Sigma_q$ we have
$G^o=IR_{\sigma}P_{\sigma}$, we get 
\begin{equation}\label{dec}
C_{q}(X,{\mathcal C} )=\bigoplus_{\sigma\in \Sigma_q}\bigoplus_{w\in
R_{\sigma}}\bigoplus_{\tau \in \Sigma [w\sigma ]}{\mathcal V}_{\tau}\ ,
\end{equation}

where $\Sigma [w\sigma ]$ is the finite set
$$
\{kw\sigma\ ; \ k\in I\} = \{ kw\sigma\ ; \ k\in I/I\cap P_{w\sigma}\} \ .
$$
 The space $\displaystyle C_{q}(X, {\mathcal C} )_{\sigma ,w}:=\bigoplus_{\tau \in \Sigma
[w\sigma ]}{\mathcal V}_{\tau}$ is a sub-$I$-module of $C_{q}(X ,{\mathcal C} )$.  An
element $(v_{\tau})_{\tau \in \Sigma [w\sigma ]}$ lies in $C_{q}(X,
{\mathcal C} )_{\sigma ,w}^{I}$ if and only if the following conditions hold:

i) $v_{kw\sigma}=kv_{w\sigma}$, $k\in I$;

ii) $v_{w\sigma}\in {\mathcal V}[w\sigma ]^{I\cap P_{w\sigma}}$.

Moreover we have:
$$
{\mathcal V}_{w\sigma}^{I\cap P_{w\sigma}}  =\left(
{\mathcal V}^{U_{w\sigma}}\right)^{I\cap P_{w\sigma}}
$$
$$
={\mathcal V}^{\langle U_{w\sigma},I\cap P_{w\sigma}\rangle}
=w{\mathcal V}^{\langle U_{\sigma}, I^{w^{-1}}\cap
P_{\sigma}\rangle}\ .
$$
According to lemma 6, we have $\langle U_{\sigma}, I^{w^{-1}}\cap
P_{\sigma}\rangle=I_1$, the Iwahori subgroup fixing $C_1
=C_{1}(\sigma ,w)$. It follows that
$$
{\mathcal V}_{w\sigma}^{I\cap P_{w\sigma}}=w{\mathcal V}^{I_1}={\mathcal V}^{I_{1}^{w}}\ .
$$
In the following we set $I_{1}^{w}=I_{\sigma ,w}$; this is the Iwahori
subgroup fixing the chamber $wC_{1}(\sigma ,w)=C(\sigma ,w)$. We have
proved:

\begin{lem} The space $C_{q}(X ,{\mathcal C} )_{\sigma ,w}^{I}$ is
given by
$$
\{(v_{\tau})_{\tau \in \Sigma [w\sigma ]}\ ; \
v_{kw\sigma}=kv_{w\sigma} , \ v_{w\sigma}\in {\mathcal V}^{I_{\sigma ,w}}\}\ .
$$
In particular we have a canonical isomorphism of vector spaces
$$
\psi_{\sigma ,w}~: \ C_{q}(X,{\mathcal C} )_{\sigma ,w}^{I}\rightarrow {\mathcal V}^{I_{\sigma ,w}}
$$
$$
(v_{\tau})_{\tau \in \Sigma [w\sigma ]}\mapsto v_{w\sigma}\ .
$$
\end{lem}

\begin{prop} The map $\varphi_{\sigma ,w}$~:
$C_{q}(X,{\mathcal C} )_{\sigma ,w}^{I}\rightarrow {\mathcal V}^I$ given by 
$$
(v_{\tau})_{\tau \in \Sigma [w\sigma ]}\mapsto \sum_{\tau \in \Sigma
[w\sigma ]}v_{\tau}
$$
is an isomorphism of vector spaces.
\end{prop}

{\it Proof}. We show that 
$$
\varphi =\varphi_{\sigma ,w}\circ \psi_{\sigma ,w}^{-1}~:\ 
{\mathcal V}^{I_{\sigma ,w}}\rightarrow {\mathcal V}^{I}
$$
is an isomorphism of vector spaces. For $v\in {\mathcal V}^{I_{\sigma ,w}}$, we
have 
$$
\varphi (v) =\sum_{k\in I/I\cap P_{w\sigma}}kv =\vert I/I\cap
P_{w\sigma}\vert e_{I}.v\ .
$$
So $\varphi =  \vert I/I\cap P_{w\sigma}\vert e_{I}$ and our
proposition follows from lemma 8.
\smallskip

 We may consider the space of $q$-chains $C_{q}({\mathcal A},{\mathcal V}^I)$ of
 the simplicial complex $\mathcal A$ with constant coefficients in ${\mathcal V}^I$.

\begin{cor} We have a natural bijection
$$
\varphi_{q}~:\ C_{q}(X ,{\mathcal C} )^I\rightarrow C_{q}({\mathcal A},{\mathcal V}^I ) 
$$
given as follows. It maps an element $(u_{\sigma})_{\sigma\in X_q}$ to
the element $(v_{\sigma})_{\sigma \in {\mathcal A}_q}$ given by 
$$
v_{\sigma}=\sum_{k\in I/I\cap P_{\sigma}}u_{k\sigma}= 
 \sum_{k\in I/I\cap P_{\sigma}} ku_{\sigma}\ , \ \sigma \in {\mathcal
 A}_q\ .
$$
\end{cor}

{\it Proof}. This is a straightforward consequence of the
decomposition (\ref{dec}) and proposition 12.

\section{Comparison of chain complexes}

 Consider the chain complex of $\mathcal A$ with (constant) coefficients in
 ${\mathcal V}^I$:

\begin{equation}\label{e}
C_{N-1} ({\mathcal A},{\mathcal V}^I ) \ {\buildrel \partial_{\mathcal A}
\over\longrightarrow}\ \cdots \ {\buildrel \partial_{\mathcal A}
\over\longrightarrow}\ C_{0}({\mathcal A},{\mathcal V}^I )\ {\buildrel
\epsilon_{\mathcal A} \over\longrightarrow}\ {\mathcal V}^I
\end{equation}

This complex is acyclic since the geometric realization of $\mathcal A$ is
contractible. In the next proposition we are going to prove that the
complexes (\ref{e}) and (\ref{d}) are isomorphic as complexes of vector
spaces; proposition 3 will then follow.
\smallskip

\begin{prop}  The following diagram is commutative:

$$
\begin{CD}
 0  & @>>> &  C_{N-1}(X, {\mathcal C} )^I & @>\Phi (\partial )>> & \cdots &
  @>\Phi (\partial )>>   &  C_{0}(X, {\mathcal C} )^I & 
@>\Phi (\epsilon )>> & {\mathcal V}^I\\
   &   &    &  & @V\varphi_{N-1}VV   &  &  &   &   & &
  @V\varphi_{0}VV  &   &
  @VV{\rm id}V    \\
0 & @>>>  & C_{N-1} ({\mathcal A},{\mathcal V}^I ) &
  @>\partial_{\mathcal A}>>  &  \cdots & @>\partial_{\mathcal A}>>
 &  C_{0}({\mathcal A},{\mathcal V}^I ) & @>\epsilon_{\mathcal A}>> &  {\mathcal V}^I
\end{CD}
$$
\end{prop}

{\it Proof}.  The proof of the proposition follows from the two
following lemmas.
\smallskip

\begin{lem} The following diagram is commutative:
$$
\begin{CD}
  C_{0}(X,  {\mathcal C} )^I & 
@>\Phi (\epsilon )>>  & {\mathcal V}^I \\
@V\varphi_{0}VV &    & @V{\rm id}VV \\
 C_{0}({\mathcal A},{\mathcal V}^I ) & @>\epsilon_{\mathcal A}>> 
&  {\mathcal V}^I 
\end{CD}
$$
\end{lem}

{\it Proof}. Let ${\bf v}=(v_{\sigma})_{\sigma\in X_0}\in C_{0}(X, {\mathcal C}
)^I$. We have 
$$
{\rm id} \circ \Phi (\epsilon )({\bf v})= \sum_{\sigma\in
X_0}v_{\sigma}\ .
$$
On the other hand, ${\bf u}=\varphi_0 ({\bf v})$ is given by
$\displaystyle u_{\sigma}=\sum_{k\in I/I\cap
P_{\sigma}}v_{k\sigma}$, $\sigma \in {\mathcal A}_{0}$. So
$$
\epsilon\circ \varphi_{0}({\bf v})=\sum_{\sigma\in {\mathcal
A}_0}\sum_{k\in I/I\cap P_{\sigma}}v_{k\sigma}
$$
and our assertion follows from the equality 
$$
X_o =\{ k\sigma \ ; \ \sigma\in {\mathcal A}_0, \ k\in I/I\cap
P_{\sigma }\}\ .
$$
Indeed any $\sigma\in X_0$ has a unique $I$-conjugate lying in $\mathcal
A$.
\smallskip

\begin{lem}  Let $q\in \{ 0,\dots ,N-2\}$. The following
diagram is commutative:
$$
\begin{CD}
  C_{q+1}(X, {\mathcal C} )^I & @>\Phi(\partial )>> &  C_{q}(X,
{\mathcal C} )^I \\
@V\varphi_{q+1}VV  &  & @V\varphi_{q}VV \\
 C_{q+1} ({\mathcal A},{\mathcal V}^I ) &  @>\partial_{\mathcal A}>>
 & C_{q} ({\mathcal A},{\mathcal V}^I 
\end{CD}
$$
\end{lem}

{\it Proof}. Let ${\bf v}=(v_{\sigma})_{\sigma\in X_{q+1}}\in
C_{q+1}(X, {\mathcal C} )^{I}$. We have
$$
\Phi (\partial )({\bf v})_{\tau}=\sum_{\sigma \supset \tau\ , \
\sigma\in X_{q+1}}[\sigma :\tau ]v_{\sigma},\ \tau \in X_{q}\ ,
$$
and
$$
\varphi_{q}\circ\Phi (\partial )({\bf v})_{\tau}=\sum_{k\in I/I\cap
P_{\tau}} \sum_{\sigma\supset k\tau \ ,\ \sigma\in
X_{q+1}}[\sigma :k\tau ]v_{\sigma}\ ,\ \tau \in {\mathcal A}_q\ .
$$
On the other hand we have
$$
\varphi_{q+1}({\bf v})_{\sigma}=\sum_{k\in I/I\cap
P_{\sigma}}v_{k\sigma}\ , \ \sigma \in {\mathcal A}_{q+1}\ , 
$$
and
$$
\partial_{\mathcal A}\circ \varphi_{q+1}({\bf
v})_{\tau}=\sum_{\sigma\supset \tau\ , \ \sigma \in {\mathcal
A}_{q+1}}[\sigma :\tau ]\sum_{k\in I/I\cap P_{\sigma}}v_{k\sigma}\ .
$$
For $\tau \in {\mathcal A}_q$ and $k\in I/I\cap P_{\tau}$, any
$\sigma\in X_{q+1}$ satisfying $\sigma \supset k\tau$ satifies
$k^{-1}\sigma\supset \tau$ and the simplex $k^{-1}\sigma$ may be
uniquely written $l\sigma_o$, for some $\sigma_o\in {\mathcal A}_{q+1}$
containing $\tau$ and some $l\in I\cap P_{\tau}/I\cap
P_{\sigma_o}$. It follows that for $\tau\in {\mathcal A}_{q}$, we
have
$$
\varphi_{q}\circ \Phi (\partial )({\bf v})_{\tau}=\sum_{k\in I/I\cap
P_{\tau}} \  \sum_{\sigma_o \supset \tau\ , \sigma_o \in {\mathcal
A}_{q+1}} \    \sum_{l\in I\cap P_{\tau}/I\cap
P_{\sigma_o}} [kl\sigma_o : k\tau ]v_{kl\sigma_o}\ .
$$
Since the action of $G^o$ preserves the orientation, for $k\in  I/I\cap
P_{\tau}$ and $l\in  I\cap P_{\tau}/I\cap
P_{\sigma_o}$, we have
$$
[kl\sigma_o :k\tau ]=[l\sigma_o :\tau ]=[\sigma_o :l^{-1}\tau
]=[\sigma_o :\tau ]\ .
$$
It follows that
$$
\varphi_{q}\circ \Phi (\partial )({\bf v})_{\tau}=
\sum_{\sigma_o \supset \tau\ , \sigma_o \in {\mathcal
A}_{q+1}} [\sigma_o :\tau ]\sum_{u\in I/I\cap
P_{\sigma_o}}v_{u\sigma_o}\ , \ \tau \in {\mathcal A}_{q}\ ,
$$
as required.

% Bibliographic references with the natbib package:
% Parenthetical: \citep{Bai92} produces (Bailyn 1992).
% Textual: \citet{Bai95} produces Bailyn et al. (1995).
% An affix and part of a reference:
%   \citep[e.g.][Ch. 2]{Bar76}
%   produces (e.g. Barnes et al. 1976, Ch. 2).

\end{document}